 \documentclass[preprint,12pt]{elsarticle}

\usepackage{amssymb}
\usepackage{amsfonts}
\usepackage{amsmath}
\usepackage{aas_macros}
\usepackage{graphics}
\usepackage{graphicx}
\usepackage{lineno}

\journal{Planetary and Space Science}
\begin{document}
\begin{frontmatter}



\title{Trajectories  and Stability  Regions  of the Lagrangian   Point  $L_1$  in the Generalized Chermnykh-Like Problem}


\author{Badam Singh Kushvah}
\ead{bskush@gmail.com,kushvah.bs.am@ismdhanbad.ac.in}
\ead[http://www.ismdhanbad.ac.in/depart/math/faculty1.htm#BSK]{http://www.ismdhanbad.ac.in/depart/math/faculty1.htm}
\address{Department of Applied Mathematics, Indian School of Mines Dhanbad-826009,INDIA,Phone:+91-326-2235765,Fax:+91-326-2296563}

\begin{abstract}
The  Lagrange point $L_1$  for  the Sun-Earth system   is considered
due to its  special  importance for the scientific community for the
design of space missions. The location of the Lagrangian points with
the trajectories and stability regions  of $L_1$ are   computed
numerically for the initial conditions very close to the point.  The
influence of belt, effect of radiation pressure due to Sun and
oblateness effect of second primary(finite body Earth) is presented
for various values of parameters. The collinear point $L_1$ is
asymptotically stable  within  a specific interval of time $t$
correspond to the values of  parameters and initial conditions.
\end{abstract}

\begin{keyword}
trajectory \sep stability \sep equilibrium points  \sep radiation
pressure \sep oblateness   \sep rtbp.
\MSC[2008] code \sep  70F15
\end{keyword}

\end{frontmatter}

\section{Introduction}
\label{Intro} The circular restricted three body problem is
modification of the three body problem where the third body is
assumed to have very small mass which is infinitesimal in comparison
to other two finite masses are  called primaries.  The restricted
three body problem is generalized to include  radiation pressure,
oblateness of the second primary and influence of the belt. Further
the  primary bodies are moving in circular orbits about their center
of mass. The well-known five equilibrium points(Lagrangian points)
that appear in the planar restricted three-body problem are very
important for astronautical applications. The collinear points are
unstable and the triangular points are stable \citet{Szebehely1967}.
In the Sun-Jupiter system several thousand asteroids, collectively
referred to as Trojan asteroids, are in orbits of triangular
equilibrium points. But collinear equilibrium points are also made
linearly stable by continuous corrections of their orbits(\lq\lq
halo orbits\rq\rq). In other words the collinear equilibrium points
are metastable points in the sense that, like a ball sitting on top
of a hill. However, in practice these Lagrange points have proven to
be very useful indeed since a spacecraft can be made to execute a
small orbit about one of these Lagrange points with a very small
expenditure of energy \citet{Farquhar1967JSpRo,Farquhar1969AsAer}.
Because of the its unobstructed view of the Sun, the Sun-Earth $L_1$
is a good place to put instruments  for doing solar science. NASA's
Genesis Discovery Mission has been there, designed completely using
invariant manifolds and other tools form dynamical systems theory.
In 1972, the International Sun-Earth Explorer (ISEE) was established
, joint project of NASA and the European Space Agency(ESA). The
ISEE-3 was launched into a halo orbit around the Sun-Earth $L_1$
point in 1978, allowing it to collect data on solar wind conditions
upstream from the Earth \citet{Farquhar1985JAnSc}. In the mid-1980s
the Solar and Heliospheric Observatory (SOHO)
\citet{Domingo1995SoPh} is places in a halo orbit around the
Sun-Earth $L_1$ position, about a million miles the Sun ward from
the Earth. They have provided useful places to \lq\lq park\rq\rq a
spacecraft for observations.

The Chermnykh's problem is a new kind of restricted three body
problem which was first time studied by \citet{Chermnykh1987}. This
problem generalizes two classical problems of Celestial mechanics:
the two fixed center problem and the  restricted three body problem.
This gives wide perspectives for applications of the problem  in
celestial mechanics and astronomy. The importance of the problem in
astronomy has been addressed  by \citet{Jiang2004IJBC}. Some
planetary systems are claimed to have discs of dust and they are
regarded to be young analogues of the Kuiper Belt in our Solar
System. If these discs are massive enough, they should play
important roles in the origin of planets\rq orbital elements. Since
the belt of planetesimal often exists within a planetary system and
provides the possible mechanism of orbital circularization, it is
important to understand the solutions of dynamical systems with the
planet-belt interaction. Chermnykh's problem  has been  studied by
many scientists such as \citet{Papadakis2005Ap&SS},
\citet{Jiang2006Ap&SS,YehJiang2006Ap&SSII},
\citet{Papadakis2007Ap&SS} and reference their  in.

The goal of present paper is to investigate the nature of collinear
equilibrium point $L_1$ because of  the interested point  to the
mission design. Although there are two new equilibrium points due to
mass of the belt(larger than 0.15)
\citet{JiangYeh2006Ap&SSI,YehJiang2006Ap&SSII} but they are left to
examine. All the results are computed numerically with the help of
computer because pure analytical methods are not suitable. The
actual trajectories and the stability  regions    of $L_1$ however
is  more complicated than the   discussed here.  But for specific
the  time intervals, and initial values,  these results provide new
information on the behavior of trajectories around the  Lagrangian
point $L_1$ for  different possible set values of the parameters.


\section{Location of Lagrangian Points}
 \label{EqJaco}
It is  supposed that the motion of an infinitesimal mass particle is
influenced by the gravitational force from primaries and a belt of
mass   $M_b$.  The  units of the mass, the distance and the time are
taken such that sum of the masses and the distance between primaries
are unities, the unit of the time i.e. the time period of $m_1$
about $m_2$ consists of $2\pi$ units such that  the Gaussian
constant of gravitational $\mathbf{ k}^{2}=1$. Then perturbed mean
motion $n$ of the primaries is given by
$n^{2}=1+\frac{3A_{2}}{2}+\frac{2M_b
r_c}{\left(r_c^2+T^2\right)^{3/2}}$, where
$T=\mathbf{a}+\mathbf{b}$, $\mathbf{a,b}$ are flatness and  core
parameters respectively which determine the density profile of the
belt,   $r_c^2=(1-\mu)q_1^{2/3}+\mu^2$,
$A_{2}=\frac{r^{2}_{e}-r^{2}_{p}}{5r^{2}}$ is the oblateness
coefficient of $m_{2}$; $r_{e}$, $r_{p}$  are the equatorial and
polar radii of $m_{2}$ respectively,  $r=\sqrt{x^2+y^2}$ is  the
distance between primaries and $x=f_1(t), y=f_2(t)$ are the
functions of the time $t$ i.e. $t$ is only independent variable. The
mass parameter is $\mu=\frac{m_{2}}{m_{1}+m_{2}}$($9.537 \times
10^{-4}$ for the Sun-Jupiter and $3.00348 \times 10^{-6}$ for the
Sun-Earth mass distributions respectively ), $q_1=1-\frac{F_p}{F_g}$
is a mass reduction factor and   $F_{p}$ is the solar radiation
pressure force which is exactly apposite to the gravitational
attraction force $F_g$.  The coordinates of $m_1$, $m_2$ are
$(-\mu,0)$, $(1-\mu,0)$ respectively. In the above mentioned
reference system and \citet*{MiyamotoNagai1975PASJ} model, the
equations of motion of the infinitesimal mass particle in the $x
y$-plane formulated as[please see
\citet{Kushvah2008Ap&SS318,Kushvah2009Ap&SS}]:
\begin{eqnarray}
\ddot{x}-2n\dot{y}&=&\Omega_x ,\label{eq:Omegax}\\
\ddot{y}+2n\dot{x}&=&\Omega_y,\label{eq:Omegay}
 \end{eqnarray}
where
\begin{eqnarray*}
\Omega_x&=& n^{2}x-\frac{(1-\mu)q_1(x+\mu)}{r^3_1}-\frac{\mu(x+\mu-1)}{r^3_2}-\frac{3}{2}\frac{\mu{A_2}(x+\mu-1)}{r^5_2}\nonumber\\
&&-\frac{M_b x}{\left(r^2+T^2\right)^{3/2}}\label{eq:Omegax1}\\
\Omega_y&=&n^{2}y
-\frac{(1-\mu)q_{1}{y}}{r^3_1}
-\frac{\mu{y}}{r^3_2}-\frac{3}{2}\frac{\mu{A_2}y}{r^5_2}\nonumber\\&&-\frac{M_b y}{\left(r^2+T^2\right)^{3/2}}\label{eq:Omegay1}\end{eqnarray*}
\begin{eqnarray}
\Omega&=&\frac{n^2(x^2+y^2)}{2}+\frac{(1-\mu)q_1}{r_1}+\frac{\mu}{r_2}+\frac{\mu
 A_2}{2r_2^3}+\frac{M_b}{\left(r^2+T^2\right)^{1/2}}\label{Omega}\\
r_1&=&\sqrt{(x+\mu)^2+y^2}, r_2=\sqrt{(x+\mu-1)^2+y^2}.\nonumber
\end{eqnarray}
From  equations (\ref{eq:Omegax}) and (\ref{eq:Omegay}),  the
Jacobian integral is given by:
\begin{eqnarray}
    E=\frac{1}{2}\left(\dot{x^2}+\dot{y}^2\right)-\Omega(x,y,\dot{x},\dot{y})=(\mbox{Constant})
\end{eqnarray}
which is related to the Jacobian constant  $C=-2E$. The location of
three collinear equilibrium points and two triangular equilibrium
points  is computed by dividing the orbital plane into three parts
$L_1,L_4(5)$: $\mu<x<(1-\mu)$, $L_2$: $(1-\mu)<x$ and $L_3$:
$x<-\mu$. For the collinear points,  an algebraic equation of the
fifth degree is solved numerically with initial  approximations to
the Taylor-series as:
\begin{eqnarray}
x(L_1)&=&1-(\frac{\mu}{3})^{1/3}+ \frac{1}{3}(\frac{\mu}{3})^{2/3}-\frac{26\mu}{27}+\dots\label{eq:xL1}\\
x(L_2)&=&1+(\frac{\mu}{3})^{1/3}+ \frac{1}{3}(\frac{\mu}{3})^{2/3}-\frac{28\mu}{27}+\dots\label{eq:xL2}\\
x(L_3)&=& -1-\frac{5\mu}{12}+ \frac{1127\mu^3}{20736}+\frac{7889\mu^4}{248832}+\dots\label{eq:xL3}\\
\end{eqnarray}
The solution of differential  equations (\ref{eq:Omegax}) and
(\ref{eq:Omegay}) is  presented as interpolation function which is
plotted  for various integration intervals by substituting specific
values of  the time $t$ and initial conditions i.e.
$x(0)=x(L_i),y(0)=0$ where $i=1--3$ and $x(0)= \frac{1}{2}-\mu,
y(0)=\pm\frac{\sqrt{3}}{2}$ for the triangular equilibrium points.
\begin{figure}
    \includegraphics[scale=.7]{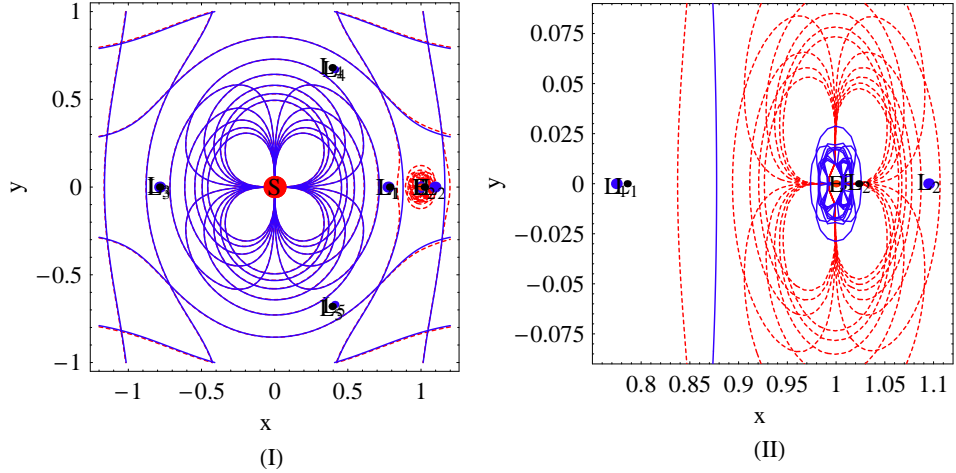}
   \caption{The position of equilibrium points  when T=0.01, $q_1=0.75$, $A_2=0.05$ and $M_b=0.4$, panel (I):Red doted  curves and blue points for Sun-Jupiter  mass distribution, blue  curves and black points for Sun-Earth mass distribution,  (II): Position of $L_1, L_2$ with respect to Jupiter's and Earth's is shown in zoom}
   \label{fig:lpoints}
\end{figure}
The equilibrium points are shown in figure \ref{fig:lpoints} in
which two panels i.e. (I)   red solid curves and blue points
correspond to the Sun-Jupiter  mass distribution and  blue dashed
curves and black  points correspond to the Sun-Earth mass
distribution. Panel (II) show the zoom of the neighborhood of $L_1,
L_2$. The numerical values of these points  are presented in Table
\ref{tab:lpts}. One can see that  the positions of $L_1, L_3$
appeared rightward and the positions of  $L_2, L_4$ ( $L_4$ is
shifted downward also) are shifted leftward in the Sun-Earth system
with respect to the position in the Sun-Jupiter system. The nature
of the $L_5$ is similar to the $L_4$. The detail behavior of the
$L_1$ with stability regions is discussed in sections.
\ref{sec:TrjL1} \& \ref{sec:stbL1}.
\begin{center}
\begin{table}
\caption{Location of equilibrium points when T=0.01, $q_1=.75$, $A_2=.05$ and $M_b=0.4$ }
\begin{tabular}{|c l l| l l|}
\hline\hline
Sun-Jupiter& &    & Sun-Earth & \\
$L_i$ & x  & y   & x  & y\\
  \hline\hline
$L_1$ &0.774577& 0  &0.78569& 0 \\ $L_2$&1.09493&0 &1.0232&0  \\ $L_3$&-0.786195&0
  & -0.785732 & 0 \\ $L_4$ &0.410603& 0.669308 & 0.393072&0.680342 \\ \hline
\end{tabular}\label{tab:lpts}
\end{table}
\end{center}
\section{Trajectory of $L_1$}
\label{sec:TrjL1} 
\begin{figure}
    \includegraphics[scale=.7]{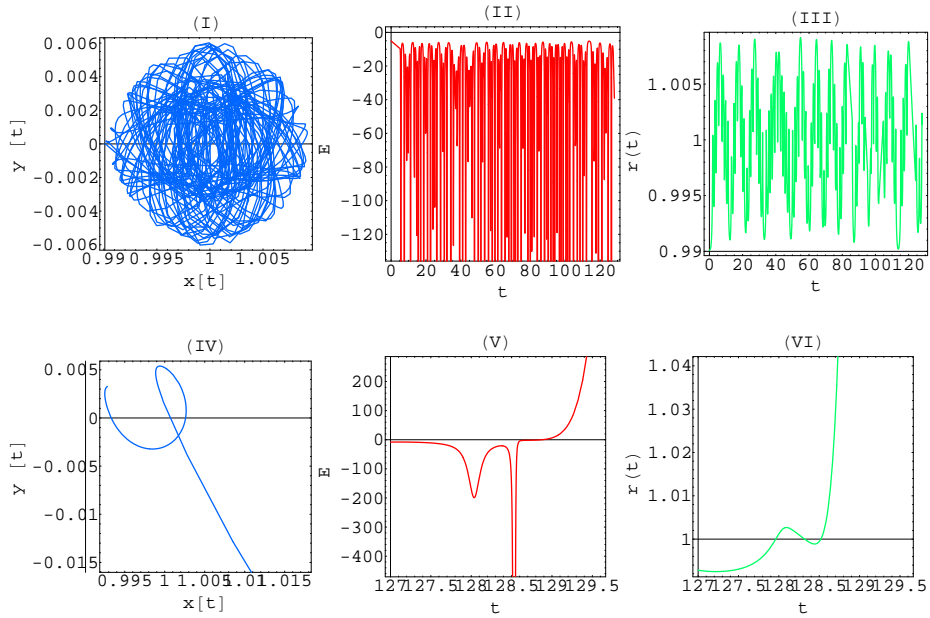}
   \caption{The Panels  (I-III):$0<t<128.23$ and (IV-VI):$127<t<129.6$ in which (I and II) show the  trajectory of $L_1$, (II and V) show  energy-versus time and (III-VI) show the local distance of trajectory at time $t$ form the initial  points i.e. $t=0$ the other parameters are T=0.01, $q_1=1$, $A_2=0$ and $M_b=0$.}
   \label{fig:trjq1a0m0}
\end{figure}

The equations (\ref{eq:Omegax}-\ref{eq:Omegay})
with initial conditions $x(0)=x(L_1), y(0)=0, x'(0)=y'(0)=0$ are
used to determine the trajectories of $L_1$ for different possible
cases. The origin of coordinate axes  is supposed to the equilibrium
point at time $t=0$ to draw the figures which show the trajectories
of the point in consideration. They are  shown in figure
\ref{fig:trjq1a0m0} with six panels i.e the panels (I-III) show the
trajectory moves about the origin ($L_1$ at $t=0$) with $x \in
(0.990093,1.00916)$ , $y\in(-0.0061448,0.00587171)$, the  energy
$E\in (-12706.5(t=22.66),-5.08226(t=0))$ and the distance
$r(t)\in(0.990093(t=0),1.00916(t=55))$. The panels
(I-III:$127<t<129.6$) show  the trajectory moves away from the
origin ($L_1$ at $t=0$) after a certain value of  the time $t$, with
$x \in (0.990093,1.00916)$ , $y\in(-0.0061448,0.00587171)$, minimum
energy $E =-1447$ found  at the time $t=128.52$, and the energy
$E>0$ for $t>128.88$.
\begin{figure}
    \includegraphics[scale=.6]{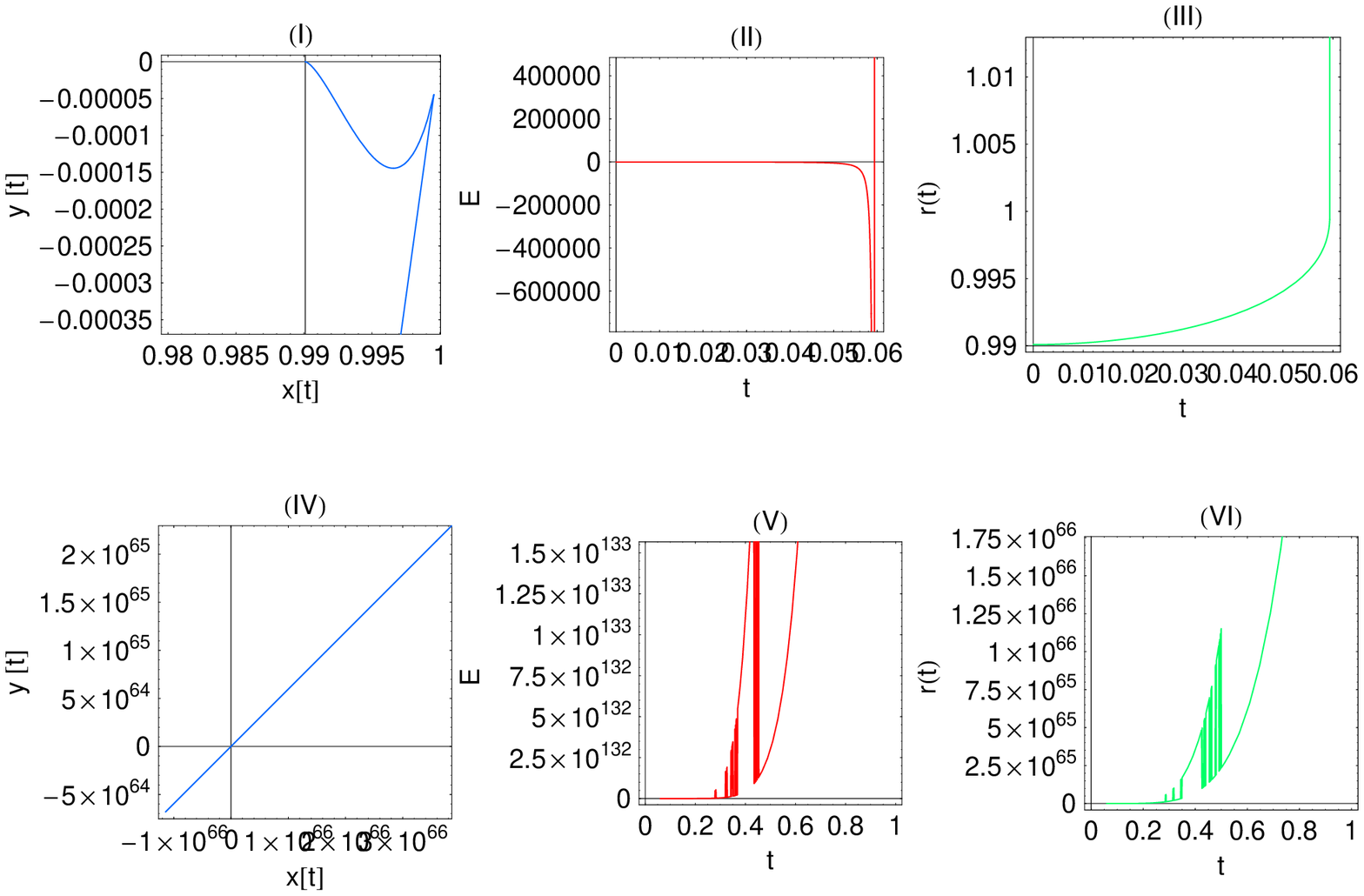}
   \caption{The Panels  (I-III):$0<t<0.06$ and (IV-VI):$.06<t<1$ in which (I and II) show the  trajectory of $L_1$, (II and V) show  energy-versus time and (III-VI) show the local distance of trajectory at time $t$ form the initial  points i.e. $t=0$ the other parameters are T=0.01, $q_1=1$, $A_2=0.05$ and $M_b=0$.}
   \label{fig:trjq1a05m0}
\end{figure}

Figure \ref{fig:trjq1a05m0} is plotted for $q_1=1, M_b=0$ and
$A_2=0.05$ with six panels (I-III:$0\leq t\leq 0.06$) and
(IV-VI:$0.06\leq t\leq 1$) which describe the effect of oblateness
of Earth to the trajectory of $L_1$. The graphs  plotted against
time which  describe  behavior of trajectories  to  equilibrium
points not the point itself is moving with time. $x=-1.91954\times
10^{-48}(t=0.06)$ to  $x=0.99405(t=.05)$  coordinate $y$ is
deceasing function of time that reach maxima -0.0000530614, at
$t=0.04$ and minima $-0.000105662$  at time $t=0.05$  again it
deceases and reach at  value $-2.5677\times 10^{47}(t=0.06)$ .
Initially energy has negative values for time $0\leq<t<0.059$
decreases with time $t$ which attains minimum value -$2.64032 \times
10^{6}(t=0.059)$ then  strictly increasing function that   attains
positive values after time $t=0.0594$. In time  interval (0.2, 0.6)
energy one time returns down that again it tend to very large
(infinite) positive value.  It is clear from panels (IV-VI) the
trajectory move far from the Lagrangian point $L_1$ after time
$t=0.0594$. The distance $r(t)$ from this point to the trajectory is
increasingly periodic for time $0<t<0.6$ then  tend to very large.

The effect of radiation pressure, oblateness and mass of the belt is
considered in figure 4, panels (I\&III) describe the trajectory and
panels (II\&IV) show the energy with respect to the time $t$. The
mass reduction factor $q_1 = 0.75$ and $M_b = 0.2$ are taken to plot
the graphs  variation in values of these parameters have similar
effect. In the panels, solid blue lines represent $A_2 = 0.25$, red
dashed lines correspond to $A_2 = 0.50$ and doted black lines for
$A_2 = 0.75$. One cane see that the trajectory move very far from
the $L_1$ the energy is positive after a certain value of the  time
$t$. Details of trajectory and energy is presented in Table
\ref{tab:trajectoryL1} for various values of parameters. One can see
that $x$ is an increasing function of the time but $y$ is an
initially decreasing function for certain values of the time, then
it becomes a strictly increasing. Similarly the energy $E$ is
negative and went downward but after some specific the time for each
cases it becomes positive and strictly increasing and attains very
large positive value.
\begin{figure}
    \includegraphics[scale=.7]{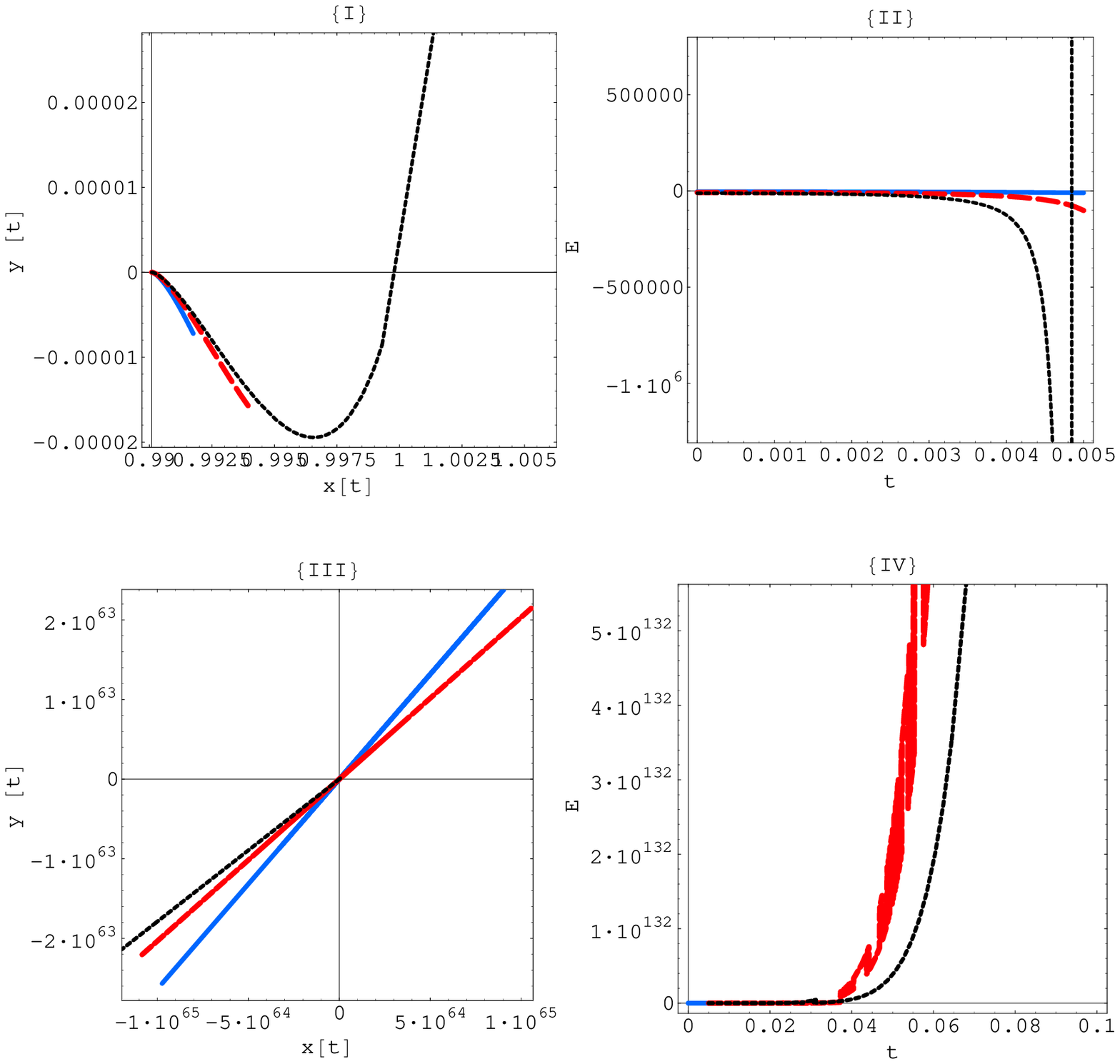}
   \caption{The Panels  (I-III):$0<t<0.005$ and (II-IV):$.00<t<1$ in which (I and II) show the  trajectory of $L_1$, (II and V) show  energy-versus time and (III-VI) show the local distance of trajectory at time $t$ form the initial  points i.e. $t=0$ the other parameters are T=0.01, $q_1=1$, $A_2=0.05$ and $M_b=0$}
   \label{fig:trjq75m2at}
\end{figure}

\begin{table}
\caption{Trajectory of $L_1$ and the energy when $T=0.01$,
$q_1=.75$,
 $M_b=0.2$ }
\begin{tabular}{|c l l  l l|}
\hline\hline $A_2$& time $t$ &  $x$  & $y$ & Energy $E$\\
  \hline\hline
 0.25 &0.000& 0.990093 & $-1.27593\times10^{-32}$& -3946.49 \\
 &0.002&0.990333   & $-4.4421 \times10^{-7}$ & -4460.93 \\
 & 0.004&0.991106&$-3.63139 \times10^{-6}$ & -6764.39\\
 &0.006&0.992646 &$-1.26827\times10^{-5}$& -17501.7 \\
  &0.008&0.996301  & $-285308\times10^{-5}$ & -544626. \\
& 0.010&$4.03799\times10^{55}$&$-1.06467\times10^{54}$
&$7.84745\times10^{117}$\\&&&&\\
 0.50 &0.000& 0.990093 & $7.51113\times10^{-33}$& -7887.36 \\
 &0.002&0.99058   & $-9.7398 \times10^{-7}$ & -10146.9\\
 & 0.004&0.992298&$-8.12979\times10^{-6}$ & -27763.1\\
 &0.006&$-8.35666\times10^{48}$ &$-1.70025\times10^{47}$& $2.17394\times10^{107}$ \\
  &0.008&$-3.05605 \times10^{56}$ & $-6.21786\times10^{54}$ & $2.74317 \times10^{119}$\\
& 0.010&$-9.05111\times10^{57}$&$-1.84154\times10^{56}$
&$6.20464\times10^{121}$\\&&&&\\
 0.75 &0.000& 0.990093 & $-1.80556\times10^{-34}$& -11828.2\\
 &0.002&0.990836  & $-1.58382 \times10^{-6}$ & -17463.6\\
 & 0.004&0.993821&$-1.33765\times10^{-5}$ & -125336.\\
 &0.006&$6.67303\times10^{55}$ &$1.19389\times10^{54}$& $4.20091\times10^{118}$ \\
  &0.008&$1.02535\times10^{58}$ & $1.83448\times10^{56}$ & $1.32379\times10^{122}$\\
& 0.010&$1.19678\times10^{59}$&$2.1412\times10^{57}$
&$6.74915\times10^{123}$\\\hline
\end{tabular}\label{tab:trajectoryL1}
\end{table}
\section{Stability of $L_1$}
\label{sec:stbL1} 
\begin{figure}
    \includegraphics[scale=.7]{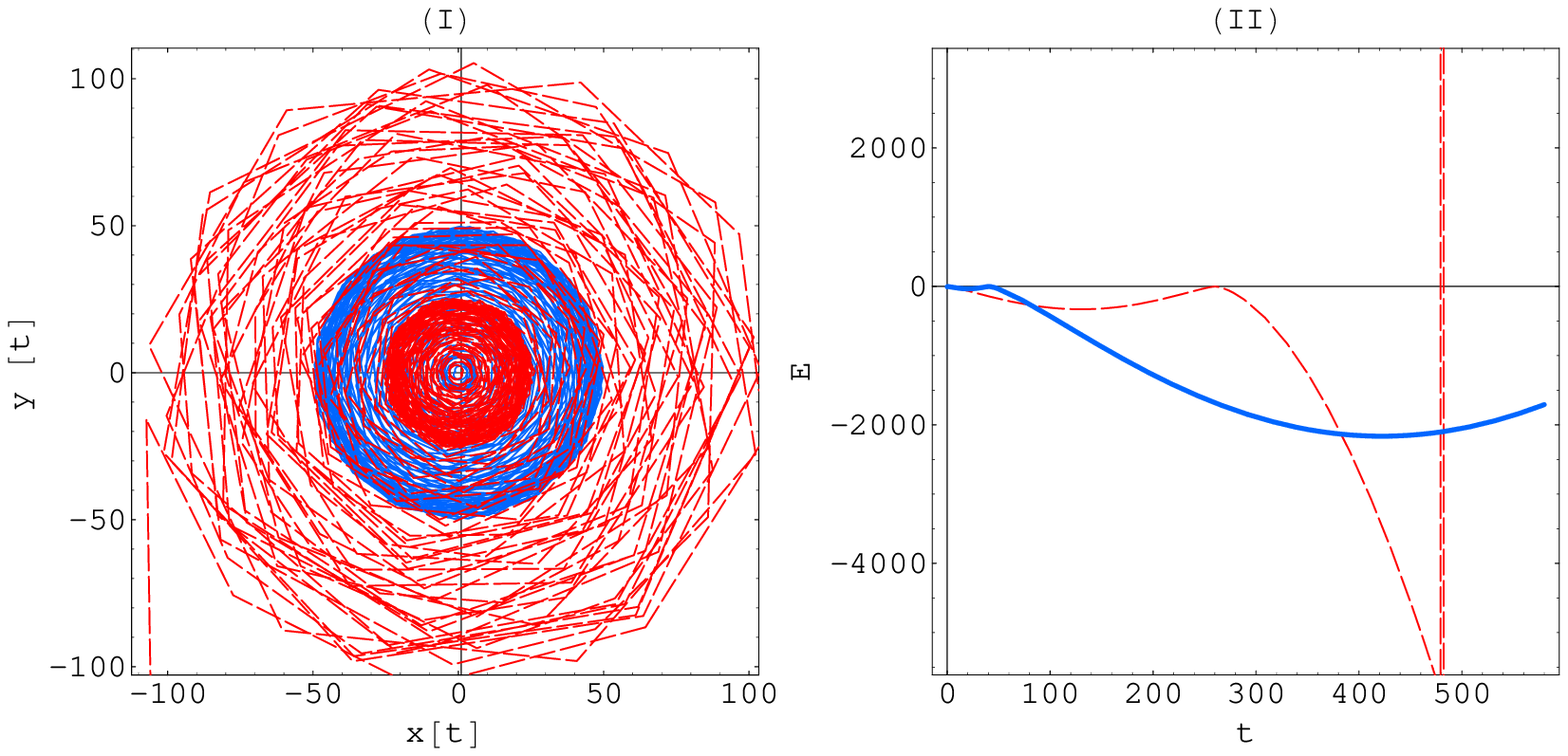}\\\includegraphics[scale=.7]{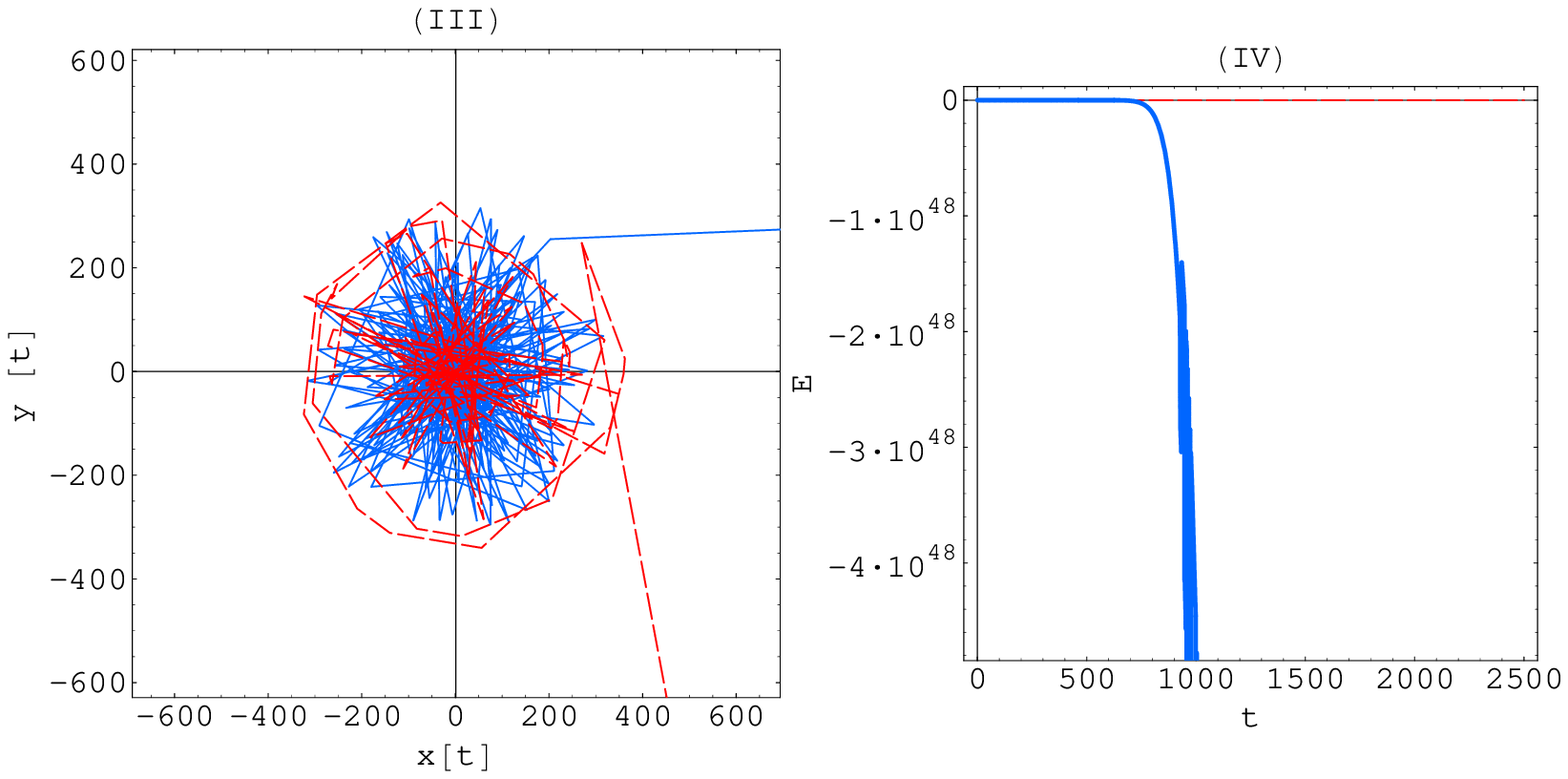}
   \caption{Show the  stability of $L_1$ with panels  (I-II):$0\leq t\lq 491,q_1=0.75,A_2=0.0$ and (III-IV):$0<t<2500, q_1=0.50, A_2=0.0$ in which blue solid curves for $M_b=0.25$, red curves for $M_b=0.50$}
   \label{fig:stb}
\end{figure}

Suppose the coordinates $(x_1, y_1)$ of $L_1$ are
initially perturbed by changing \(x(0) = x_1+\epsilon \cos(\phi),
y(0) = y_1+\epsilon\sin(\phi)\) where \( \phi
=\arctan\left(\frac{y(0)-y_1}{x(0)-x_1}\right)\in  (0, 2\pi), 0 \leq
\epsilon= \sqrt{(x(0) - x_1)^2 + (y(0) - y_1)^2} < 1\). The  $\phi$
indicates the direction of the initial position vector in the local
frame. If the $\epsilon = 0$ means there is no perturbation. It is
supposed that the $\epsilon = 0.001$ and the $\phi= \frac{\pi}{4}$
to examine the stability of $L_1$. Figure \ref{fig:stb} show the
path of test particle and its energy with four panels i.e. the
panels (I\&III):$ q_1 = 0.75, 0.50,A_2 = 0.0$, in (I) trajectory of
perturbed $L_1$ moves in chaotic-circular path around initial
position without deviating far from it, then steadily move out of
the region. In (III) the test particle move in stability region and
returns repeatedly on its initial position. The blue solid curves
represent $M_b = 0.25$ and dashed curves represent $M_b = 0.50$. It
is clear form panel (III) that bounded region for $M_b = 0.25$ is $t
< 2500$ and for $M_b = 0.50, t < 2600$.

\begin{figure}
    \includegraphics[scale=.7]{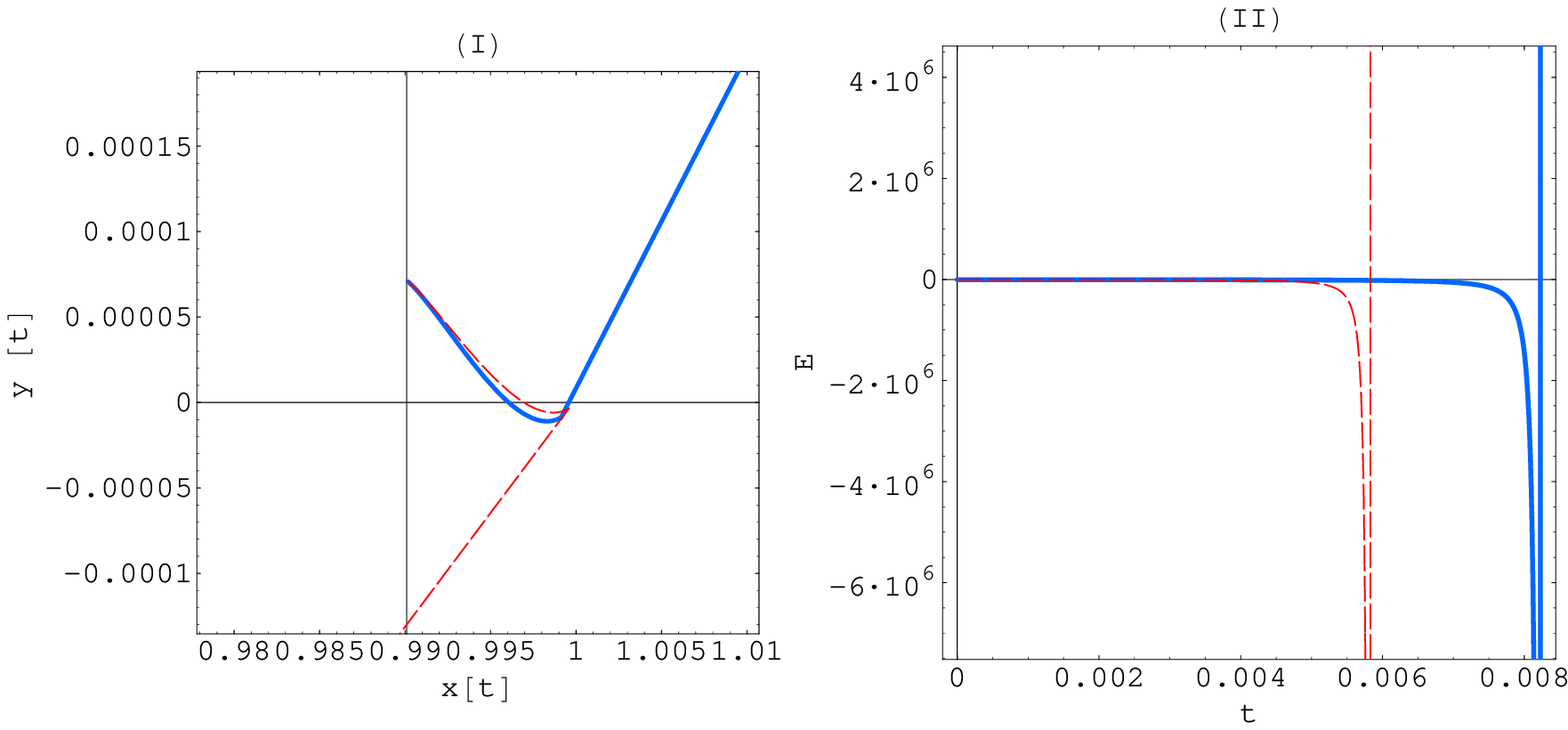}\\
   \caption{Effect of oblateness coefficient $A_2$ on the stability of $L_1$ panel (I) trajectory  (II) energy of perturbed  point $L_1$  in which blue solid curves for $A_2=0.25$, red curves for $A_2=0.50$}
   \label{fig:stbOBLT}
\end{figure}

The effect of oblateness of the second primary is shown in figure
\ref{fig:stbOBLT} when $q_1=0.75 , M_b=0.25$. The panel (I) shows
the  trajectory of perturbed point $L_1$ and (II) shows the energy
of that point. The blue doted lines correspond to $A_1=0.25$ and red
lines for $A_2=0.50$. One can see that the oblate effect is very
powerful on the trajectory and stability of $L_1$. When $A_2=0.0$
the $L_1$ is asymptotically stable  for the value of $t$ which lies
within a certain interval. But if oblate effect of second primary is
present($A_2\neq 0$), the stability region of $L_1$ disappears when
this effect is increases. Further all the results presented in the
manuscripts  are similar to the results  obtained by
\cite{Grebennikov2007CMMPh}, \citet{Kushvah2009RAA}.

\section{Conclusion}
The numerical computation presented in the manuscript  provides
remarkable results to design trajectories  of Lagrangian point $L_1$
which helps us to make comments on the stability(asymptotically) of
the point. We obtained the intervals of the  time where trajectory
continuously  moves around the $L_1$,  does not deviate far from the
point but tend to approach (for some cases) it, the energy of
perturbed point is negative for these intervals, so we conclude that
the point is asymptotical stable. More over we have seen that after
the specific time intervals the trajectory of perturbed point depart
from the neighborhood and goes away from it, in this case the energy
also  becomes positive, so the Lagrangian point  $L_1$ is unstable.
Further the trajectories and the stability regions are affected by
the radiation pressure, the oblateness of the second primary and
mass of the belt.
\section{Acknowledgements}
I  wish to express my  sincere thanks to the Director, Indian School
of Mines Dhanbad, the Head of the Department of  Applied Mathematics
for providing all accessary research facilities in the Department. I
also thanks to  Professor B. Ishwar, B.R.A. Bihar University
Muzaffarpur, all my colleagues and special thanks to Dr. P.S. Rao
for the valuable comments and suggestions during the preparation of
this article.


\bibliographystyle{model1a-num-names}





\end{document}